\documentclass[11pt]{article}
\usepackage{amsfonts,amssymb}
\usepackage{amsmath}
\usepackage{amsthm}
\usepackage[matrix,arrow,curve]{xy}
\newtheorem{theorem}{Theorem}
\newtheorem{lemma}{Lemma}
\newtheorem{prop}{Proposition}
\newtheorem{corollary}{Corollary}
\newtheorem{definition}{Definition}

\title{Spectrum of the Second Variation}
\author{A.~A.~Agrachev\thanks{SISSA, Trieste, MIAN, Moscow, and PSI RAS, Pereslavl. This work was supported by the Russian Sience Foundation under grant No 17-11-01387.}}
\date{}
\begin{document}
\maketitle

\begin{abstract}
Second variation of a smooth optimal control problem at a regular extremal is a symmetric Fredholm operator.
We study asymptotics of the spectrum of this operator and give an explicit expression for its determinant in terms of solutions
of the Jacobi equation. In the case of the least action principle for the harmonic oscillator we obtain a classical Euler identity
$\prod\limits_{n=1}^\infty\left(1-\frac{x^2}{(\pi n)^2}\right)=
\frac{\sin x}{x}$. General case may serve as a rich source of new nice identities.
\end{abstract}

\section{Introduction}
This research was initially motivated by the analysis of asymptotic expansions related to semigroups generated
by hypo-elliptic operators. These expansions can be often interpreted as infinite-dimensional versions of the
standard asymptotics for expressions of the form $\int\limits_{\mathbb R^N}a(u)e^{\frac{\varphi(u)}{-t}}du$ or
$\int\limits_{\mathbb R^N}a(u)e^{\frac{\varphi(u)}{it}}du$ as $t\to 0$. It is well-known that the asymptotics is localized in the critical points $x_i$ of $\varphi$, and the quantities $(\det D^2_{x_i})^{-\frac 12}$
play a crucial role in the expansions. In the infinite-dimensional version we deal with a ``path integral'', and
$\varphi$ is an ``action functional'' whose critical points are extremals of the appropriate optimal control problem. So we try to analyse and compute $\det D^2_{x_i}$ in this case. What we obtain, can be considered as
a modification and generalization of the classical ``Hill's determinant'' in Rational Mechanics (see \cite{BoTr}).

Before the calculation of the determinant, we study asymptotics of the spectrum of the second variation $D_u^2\varphi$ that is a symmetric Fredholm operator of the form $I+K$, where $K$ is a compact Hilbert--Schmidt operator. The point is
that $K$ is usually NOT a trace class operator so that the trace of $K$ and the determinant of $I+K$ are not
well-defined in the standard sense.

Anyway, a specific structure of the spectral asymptotics allows to define and compute everything. This asymptotic
is described in Theorem~1. Actually, there is an important symmetry in the asymptotics that leads to a cancellation
of slow convergent to zero positive and negative terms when we compute the sum or the product.

There is a general reason for this symmetry. It concerns an evolutionary nature of the optimal control problems:
the space of available control functions grows with time. This evolutionary structure results in the following
property (I would call it ``causality'') of the operator $K$: there exists a Volterra operator $V$ on
$L^2([0,1];\mathbb R^m)$ and a finite codimension subspace of $L^2([0,1];\mathbb R^m)$ such that
$\langle Ku,u\rangle=\langle Vu,u\rangle$ for any $u$ from the subspace. In other words, our symmetric operator imitates a triangle one on a subspace of finite codimension. We expect that all causal operators have spectral
asymptotics similar to the asymptotics from Theorem~1, although we do not study general causal operators in this
paper.

In Theorem~2 we give an explicit expression for $\det(I+K)$ in terms of solutions of the Jacobi equation that generalize classical Jacobi fields. We also give an explicit integral expression for $\mathrm{tr}\,K$;
moreover, integral formulas for all elementary symmetric functions of the eigenvalues of $K$ can be recovered from
the main determinantal formula.

I hope that this result would serve as an effective summation method, a way to get explicit expressions for infinite products and sums of interesting series. A simple example: for the 1-dimensional linear control system
$\dot x=ax+u$ with the quadratic cost $\varphi(u)=\int_0^1u^2(t)-(a^2+b^2)x^2(t)\,dt$ our determinantal
identity reads: $\prod\limits_{n=1}^\infty\left(1-\frac{a^2+b^2}{a^2+(\pi n)^2}\right)=
\frac{a\sin b}{b\,\mathrm{sh}\,a}$; the case $a=0$ corresponds to the famous Euler identity
$\prod\limits_{n=1}^\infty\left(1-\frac{b^2}{(\pi n)^2}\right)=
\frac{\sin b}{b}$.

\section{Preliminaries}

We consider a smooth control system of the form:
$$
\dot q=f^t_u(q),\quad q\in M,\ u\in U,  \eqno (1)
$$
with a fixed initial point $q_0\in M$. Here $M$ and $U$ are smooth manifolds (without border) and the vector $f^t_u(q)\in T_qM$ smoothly depends on $(q,u)\in M\times U$ and is measurable bounded with respect to $t\in[0,1]$.

Let $u(\cdot):[0,1]\to U$ be a measurable map such that $u([0,1])$ is contained in a compact subset of $U$. The Cauchy problem
$$
\dot q=f^t_{u(t)}(q),\quad q(0)=q_0,
$$
has a unique Lipschitzian solution $t\mapsto q(t;u(\cdot))$ defined on an interval in $\mathbb R$. We say that $u(\cdot)$ is an admissible control and $q(\cdot;u(\cdot))$ is a corresponding admissible trajectory if the domain
of this solution contains $[0,1]$. We denote the set of all admissible controls by $\mathcal U$; then
$\mathcal U$ is an open subset of $L^\infty\left([0,1];U\right)$. Hence $\mathcal U$ is a smooth Banach manifold modelled on the space $L^\infty\left([0,1];\mathbb R^{\dim U}\right)$.

Given $t\in[0,1]$ we define the ``evaluation map'' $F_t:\mathcal U\to M$ by the formula
$F_t(u(\cdot))=q(t;u(\cdot))$; then $F_t$ is a smooth map from the Banach manifold $\mathcal U$ to $M$.

Let $\ell:M\times U\to\mathbb R$ be a smooth ``Lagrangian''. We consider functionals
$\varphi_t:\mathcal U\to\mathbb R,\ 0\le t\le 1,$ defined by the formula:
$$
\varphi_t(u(\cdot))=\int_0^t\ell(q(\tau;u(\cdot)),u(\tau))\,d\tau.  \eqno (2)
$$

\begin{definition}
We say that $u\in\mathcal U$ is a normal extremal control if there exists $\lambda_1\in T^*_{F_1(u)}M$ such that
$\lambda_1D_uF_1=d_u\varphi_1$; here $\lambda_1D_uF_1$ is the composition of $D_uF_1:T_u\mathcal U\to T_{F_1(u)}M$
and $\lambda_1:T_{F_1(u)}M\to \mathbb R$.
We say that a normal extremal control is strictly normal if it is a regular point of $F_1$.
\end{definition}

Strictly normal extremal controls are just critical points of $\varphi\bigr|_{F^{-1}_1(q)}$, $q\in M$.
Let $t\in[0,1]$; it is easy to see that the restriction $u|_{[0,t]}$ of a normal extremal control $u$ is also a normal extremal control in the following sense:
$$
\exists\ \lambda_t\in T^*_{F_t(u)}M\quad \mathrm{such\ that}\quad \lambda_tD_uF_t=d_u\varphi_t. \eqno (3)
$$
The ``differentiation of the identity (3) with respect to $t$'' leads to the following Hamiltonian characterization of normal extremal controls.

A family of Hamiltonians $h^t_u:T^*M\to \mathbb R,\ u\in U$, is defined by the formula:
$$
h^t_u(\lambda)=\langle \lambda,f^t_u(q)\rangle-\ell(q,u),\quad g\in M,\ \lambda\in T^*_qM.
$$
Let $\sigma$ be the canonical symplectic form on $T^*M$ and $\pi:T^*M\to M$ be the standard projection,
$\pi(T^*_qM)=q$. Recall that $\sigma=ds$ where $s$ is a Liouville (or tautological) 1-form,
$\langle s_\lambda,\eta\rangle=\langle\lambda,\pi_*\eta\rangle,\ \forall\,\lambda\in T^*M,\ \eta\in T_\lambda(T^*M)$. Given a smooth function $h:T^*M\to\mathbb R$, the Hamiltonian vector field $\vec h$ on $T^*M$
is defined by the identity: $dh=\sigma(\cdot,\vec h)$. 

\begin{prop} Let $\tilde u\in\mathcal U$ and $\tilde q(t)=q(t;\tilde u),\ 0\le t\le 1$; then $\tilde u$ is a normal extremal control if and only if there exists $\tilde\lambda_t\in T^*_{\tilde q(t)}M$ such that
$$
\dot{\tilde\lambda}_t=\vec h^t_{\tilde u(t)},\quad \frac{\partial h^t_u(\tilde\lambda_t)}{\partial u}\Bigr|_{u=\tilde u(t)}=0,\quad 0\le t\le 1. \eqno (4)
$$
\end{prop}
This statement is standard, it is actually a weak version of the Pontryagin maximum principle.

A Lipschitzian curve $\tilde\lambda_t,\ 0\le t\le 1$, in $T^*M$ that satisfies relations (4) is called a normal extremal associated to the normal extremal control $\tilde u$. The time-varying Hamiltonian system
$\dot\lambda=\vec h^t_{\tilde u(t)}(\lambda)$ defines a flow
$$
\tilde\Phi^t:T^*M\to T^*M,\quad \tilde\Phi^t:\lambda(0)\mapsto\lambda(t),
$$
where $\dot\lambda(\tau)=\vec h^t_{\tilde u(\tau)}(\lambda(\tau)),\ 0\le\tau\le t$.
Obviously, $\tilde\Phi^t(\tilde\lambda_0)=\tilde\lambda_t$. Moreover, $\pi_*\vec h^t_{\tilde u(t)}(\lambda)=f^t_{\tilde u(t)}(q),\ \forall\,q\in M,\ \lambda\in T^*_qM$.

Let $\tilde P^t:M\to M$ be the flow generated by the time-varying system $\dot q=f^t_{\tilde u(t)}(q)$, i.\,e.
$\tilde P^t:q(0)\mapsto q(t)$, where $\dot q(\tau)=f^\tau_{\tilde u(\tau)}(q(\tau)),\ 0\le\tau\le t$. It follows
that $\tilde\Phi^t$ are fiberwise transformations and $\tilde\Phi^t(T_q^*M)=T^*_{\tilde P^t(q)}M$. Moreover,
the restriction of $\tilde\Phi^t$ to $T_q^*M$ is an affine map of the vector space $T^*_qM$
on the vector space $T^*_{\tilde P^t(q)}M$.

In what follows, we assume that $\tilde u$ is a normal extremal control and $\tilde\lambda_t$ is a normal extremal, they are fixed until the end of the paper. Then $\tilde u(t)$ is a critical point of the function
$u\mapsto h^t_u(\tilde\lambda_t),\ u\in U$ for any $t\in[01]$. Hence the Hessian
$\frac{\partial^2h^t_{\tilde u(t)}}{\partial u^2}(\tilde\lambda_t)$ is a well-defined quadratic form on $T_{\tilde u(t)}U$ or, in other words, a self-adjoint linear map from $T_{\tilde u(t)}U$ to $T^*_{\tilde u(t)}U$.

Now we consider Hamiltonian functions
$$
g^t_u=(h^t_u-h^t_{\tilde u(t)})\circ\tilde\Phi^t,\quad u\in U,\ t\in[0,1].
$$
A time-varying Hamiltonian vector field $\vec g^t_u$ generates the flow $(\tilde\Phi^t)^{-1}\circ\Phi^t_u$, where
$\Phi^t_u$ is the flow generated by the field $\vec h^t_u$.
We have:
$$
g^t_{\tilde u(t)}\equiv 0,\quad \frac{\partial g^t_{\tilde u(t)}}{\partial u}(\tilde\lambda_0)=0,\quad
\frac{\partial^2g^t_{\tilde u(t)}}{\partial^2 u}(\tilde\lambda_0)=\frac{\partial^2h^t_{\tilde u(t)}}{\partial^2 u} (\tilde\lambda_t).
$$
We introduce a simplified notation $H_t\doteq\frac{\partial^2h^t_{\tilde u(t)}}{\partial^2 u}(\tilde\lambda_t)$;
recall that $H_t:T_{\tilde u(t)}U\to T^*_{\tilde u(t)}U$ is a self-adjoint linear map.

More notations. Consider the map $u\mapsto \vec g^t_u(\tilde\lambda_0)$ from $U$ to $T_{\tilde\lambda_0}(T^*M)$
We denote by $Z_t$ the differential of this map at $\tilde u(t),\ Z_t\doteq \frac{\partial \vec g^t_{\tilde u(t)}}{\partial u}(\tilde\lambda_0)$; then $Z_t$ is a linear map from $T_{\tilde u(t)}U$ to $T_{\tilde\lambda_0}(T^*M)$. We also set $X_t=\pi_*Z_t$; then $X_t$ is a linear map from $T_{\tilde u(t)}U$ to
$T_{\tilde q(t)}M$. Finally, we denote by $J:T_{\tilde\lambda_0}(T^*M)\to T_{\tilde\lambda_0}^*(T^*M)$ the anti-symmetric linear map defined by the identity $\sigma_{\tilde\lambda_0}(\cdot,\cdot)=\langle J\cdot,\cdot\rangle$.

Note that $T_{u(\cdot)}\mathcal U$ is the space of measurable bounded mappings $t\mapsto v(t)\in T_{u(t)}U,\ 0\le t\le 1$. A simple rearrangement of the standard formula for the first variation in the optimal control theory (see for instance the textbook \cite{AgSa}) gives the identity:
$$
\left(D_{\tilde u}F_t\right)v=\tilde P_*^t\int_0^tX_\tau v(\tau)\,d\tau.
$$

Assume that $\tilde u$ is a regular point of $F_t$ (i.\,e. $D_{\tilde u}F_t$ is surjective); then $\tilde u$ is a critical point of
$\varphi_t|_{F_t^{-1}(\tilde q(t))}$.  We have:
$$
T_{\tilde u}(F_t^{-1}(\tilde q(t))=\ker D_{\tilde u}F_t=\left\{v\in T_{\tilde u}\mathcal U : \int_0^tX_{\tau}v(\tau)\,d\tau=0\right\}.
$$

The Hessian of $\varphi_1|_{F_1^{-1}(\tilde q(1))}$ at $\tilde u$ is a quadratic form
$D^2\varphi_1:\ker D_{\tilde u}F_1\to\mathbb R$. This is the ``second variation'' of the optimal control problem at $\tilde u$, the main object of this paper. It has the following expression (see \cite{AgSa}):
$$
D^2_{\tilde u}\varphi_1(v)=-\int_0^1\langle H_tv(t),v(t)\rangle\,dt-
\int_0^1\left\langle J\int_0^tZ_\tau v(\tau)\,d\tau,Z_tv(t)\right\rangle\,dt, \eqno (5)
$$
where $v\in T_{\tilde u}\mathcal U$ and $\int_0^1X_tv(t)\,dt=0$.

\begin{definition} The extremal $\tilde\lambda_t,\ 0\le t\le 1,$ is regular if $H_t$ is invertible for any $t\in [0,1]$ with a uniformly bounded inverse.
\end{definition}

In what follows, we assume that the reference extremal is regular. Moreover, we assume that
$H_t<0$\,\footnote{We say that a self-adjoint map is positive (negative) if the corresponding quadratic form is positive (negative).} for any $t\in[0,1]$. This last assumption is motivated by the classical Legendre condition:
for a regular extremal, integral quadratic form (5) has a finite negative inertia index if and only if
$H_t<0$, for any $t\in[0,1]$.

Finally we introduce a ``Gramm matrix'', a self-adjoint linear map $\Gamma_t:T^*_{q_0}M\to T_{q_0}M$ defined by the formula: $\Gamma_t=-\int_0^tX_\tau H^{-1}_t X_\tau^*\,d\tau$. We see that $\tilde u$ is a regular point of $F_t$ if and only if $\Gamma_t$ is invertible.
\smallskip

{\bf Remark 1.} There is an apparently more general and natural way to define a control system. Let me briefly describe it confining within the time-invariant case. Indeed, instead of
the product $M\times U$ we may consider a locally trivial bundle over $M$ with a typical fiber $U$. Then $f$ is a
smooth fiberwise map from this bundle to $TM$ that sends the fiber $U_q$ into $T_qM,\ q\in M$. In this setting,
the control and the correspondent trajectory is somehow a unique object: $u(t)\in U_{q(t)}$ and
$\dot q(t)=f(u(t))$.

The situation is reduced to what we did before this remark if the bundle is trivial.
In the general case, we simply trivialize the bundle in a neighborhood of the reference trajectory. To be more precise, we first take $\mathbb R\times M$, where $\mathbb R$ is the time axis and trivialize the pullback bundle over $\mathbb R\times M$ in a neighborhood of the graph of the reference trajectory (the graph does not have
self-intersections while the original trajectory might have). Moreover, all statements of this paper use only
$H_t$, $Z_t$, and $X_t$, $0\le t\le 1$, and these quantities depend only on the trivialization of the vector bundle  $T_{\tilde u(t)}U_{q(t)}$ along (the graph of) the curve $u(t),\ 0\le t\le 1$. A trivialization jf the vector bundle along such a curve is naturally achieved by a parallel transport in virtue of a linear connection on the vector bundle with the fibers
$T_uU_q,\ q\in M,u\in U_q$ and the base $\bigcup\limits_{q\in M}U_q$.

So, a linear connection on this vector bundle is actually all we need to write all the formulas. I leave to an interested reader to write them. I have decided not to follow this way in order to avoid an unnecessary language
complication and to make the paper affordable for a larger audience.

\section{Main results}

According to our construction, the space $T_{\tilde u}\mathcal U$ consists of the $L^\infty$-maps
$t\mapsto v(t)\in T_{\tilde u(t)}U,\ 0\le t\le 1$. At the same time, linear map
$D_{\tilde u}F_1:T_{\tilde u}\mathcal U\to T_{\tilde q(1)}M$ and quadratic
form $D^2_{\tilde u}\varphi_1$ are continuous in a weaker topology $L^2$. Let $\mathcal V$ be the closure of
$\ker D_{\tilde u}F_1$ in the topology $L^2$. Then $\mathcal V$ is a Hilbert space equipped with a Hilbert structure
$$
<v_1|v_2>\doteq\int_0^1\langle -H_tv_1(t),v_2(t)\rangle\,dt.
$$

Formula (5) implies that
$$
D^2_{\tilde u}\varphi_1(v)=<(I+K)v|v>,\quad v\in\mathcal V, \eqno (6)
$$
where $K$ is a compact symmetric operator on $\mathcal V$. In particular, the spectrum of $K$ is real, the only limiting point of the spectrum is 0, and any nonzero eigenvalue has a finite multiplicity.

What can we say about
the ``trace'' of $K$ and the ``determinant'' of $I+K$?
Let us consider a simple example, the least action principle for a charged particle in the plane in a constant magnetic field.

\smallskip
{\bf Example 1.} Let $M=U=\mathbb R^2,\ q=(q^1,q^2),\ u=(u^1,u^2),\ f_u(q)=u,\
\ell(q,u)=\frac 12|u|^2+r(q^1u^2-q^2u^1),\ q_0=0,\ \tilde u(t)=0$. Then
$$
T^*M=\mathbb R^2\times\mathbb R^2=\{(p,q): p\in\mathbb R^2,\ q\in\mathbb R^2\},\quad J(p,q)=(-q,p),
$$
$$
h_u(p,q)=\langle p,u\rangle-\frac 12|u|^2-r(q^1u^2-q^2u^1),
$$
$$
h_{\tilde u(t)}=0,\quad
\tilde\lambda_t=(\tilde p(t),\tilde q(t))=(0,0),\quad g_u^t=h_u.
$$
We have $H_tv=-v,\ Z_tv=(rv^2,-rv^1;v^1,v^2),\ X_tv=v,\ t\in[0,1],\ v\in\mathbb R^2$. Then
$$
\mathcal V=\{v\in L^2([0,1];\mathbb R^2) : \int_0^1v(t)\,dt=0\}.
$$
It is convenient to identify $\mathbb R^2$ with $\mathbb C$ as follows: $(v^1,v^2)=v^1+iv^2$. An immediate calculation gives the following expression for the operator $K$:
$$
Kv(t)=\int_0^t2riv(\tau)\,d\tau-\int_0^1\int_0^t2riv(\tau)\,d\tau dt.
$$
The eigenfunctions of this operator have a form $t\mapsto ce^{2\pi nit},\ 0\le t\le 1,\ c\in\mathbb C,\ n=\pm 1,\pm 2,\ldots$,
where the eigenfunction $ce^{2\pi nit}$ corresponds to the eigenvalue $\frac r{\pi n}$.

\smallskip
We see that even in this model example the eigenvalues of $K$ do not form an absolutely convergent series and $K$
is not a trace class operator. On the other hand, the next theorem implies that a ``principal value'' of such a series does exist at least if the data are piece-wise real analytic.

Consider the operator $H_t^{-1}Z_t^*JZ_t:T_{\tilde u(t)}U\to T_{\tilde u(t)}U$. This operator is associated to an
anti-symmetric bilinear form and has only purely imaginary eigenvalues. We denote by $\bar\zeta_t$ the sum of positive eigenvalues of the operator $iH_t^{-1}Z_t^*JZ_t$ counted according to multiplicity (here $i$ is the imaginary unit).

\smallskip
{\bf Remark 2.} Analysing notations of Section~2 we see that $\bar\zeta_t$ is the sum of positive roots (counted according to multiplicity) of the equation
$$
\det\left(\Bigl\{\frac{\partial h^t_{\tilde u(t)}}{\partial u},
\frac{\partial h^t_{\tilde u(t)}}{\partial u}\Bigr\}(\tilde\lambda_t)+si\frac{\partial^2h^t_{\tilde u(t)}}{\partial u^2}(\tilde\lambda_t)\right)=0
$$
with unknown $s$. Here $\{\cdot,\cdot\}$ is the Poisson bracket so that
$$
(v,v')\mapsto \Bigl\{\frac{\partial h^t_{\tilde u(t)}}{\partial u}v,
\frac{\partial h^t_{\tilde u(t)}}{\partial u}v'\Bigr\}(\tilde\lambda_t),\quad v,v'\in T_{\tilde u(t)}U,
$$
is an anti-symmetric bilinear form and $\frac{\partial^2h^t_{\tilde u(t)}}{\partial u^2}(\tilde\lambda_t)$, is a symmetric bilinear form on $T_{\tilde u(t)}U$.

\medskip
Let $\mathrm{Sp}(K)\subset\mathbb R$ be the spectrum of the operator $K,\ \mathrm{Sp}(K)\setminus\{0\}=\mathrm{Sp}_+(K)\cup\mathrm{Sp}_-(K)$, where $\mathrm{Sp}_\pm(K)\subset\mathbb R_\pm$.  Given $\alpha\in\mathrm{Sp}(K)\setminus{0}$, we
denote by $m_\alpha$ the multiplicity of the eigenvalue $\alpha$. Moreover, if $\mathrm{Sp}_\pm(K)$ is an infinite set, then we introduce a natural ordering of
$\mathrm{Sp}_\pm(K)$ that is a monotone decreasing sequence $\alpha_n,\ n\in\mathbb Z_\pm$, with the following properties:
$$
\bigcup\limits_{n\in\mathbb Z_\pm}\{\alpha_n\}=\mathrm{Sp}_\pm(K),\quad \#\{n\in\mathbb Z_\pm: \alpha_n=\alpha\}=m_\alpha. \eqno (7)
$$

\begin{theorem}If $\bar\zeta_t\equiv 0$, then $\alpha_n=O(|n|^{-2})$ as $n\to\pm\infty$. If $\bar\zeta_t$ is not identical zero and
 $H_t$ and $Z_t$ are piecewise real analytic with respect to $t$ then $\mathrm{Sp}_+(K)$ and $\mathrm{Sp}_-(K)$ are both infinite and
$$
\alpha_n=\frac 1{\pi n}\int_0^1\bar\zeta_t\,dt+O\left(|n|^{-5/3}\right)\quad \mathrm{as}\quad n\to\pm\infty. \eqno (8)
$$
\end{theorem}

{\bf Remark 3.} It is reasonable to expect that the statement of the theorem is valid without the piecewise-analyticity assumption. Moreover, the order $n^{-5/3}$ of the remainder term is certainly not optimal and perhaps can be substituted by $n^{-2}$. Anyway, the stated result is quite sufficient for our purposes while the proof of a stronger one would require a more sophisticated technique.

\smallskip
A cancellation of slow convergent to zero terms of the opposite sign in the expansion (8) gives the following:
\begin{corollary} The depending on $\varepsilon>0$ families of real numbers
$$
\sum\limits_{\alpha\in\mathrm{Sp}(K)\atop |\alpha|\ge\varepsilon}m_\alpha\alpha, \quad
\prod\limits_{\alpha\in\mathrm{Sp}(K)\atop |\alpha|\ge\varepsilon}(1+\alpha)^{m_\alpha}
$$
have finite limits as $\varepsilon\to 0$.
\end{corollary}
We use natural notations for these limits:
$$
\mathrm{tr}K=\lim\limits_{\varepsilon\to 0}\sum\limits_{\alpha\in\mathrm{Sp}(K)\atop |\alpha|\ge\varepsilon}m_\alpha\alpha, \quad
\det(I+K)=\lim\limits_{\varepsilon\to 0}\prod\limits_{\alpha\in\mathrm{Sp}(K)\atop |\alpha|\ge\varepsilon}(1+\alpha)^{m_\alpha}.
$$

We are going to compute these trace and determinant in terms of $H_t,\ Z_t,\ X_t$ and solutions of the following
{\it Jacobi system}:
$$
\dot\eta=-Z_\tau H^{-1}_\tau Z^*_{\tau}J\eta,\quad \eta(t)\in T_{\tilde\lambda_0}(T^*M),\ 0\le \tau\le 1. \eqno (9)
$$
This is a linear time-varying Hamiltonian system in $T_{\tilde\lambda_0}(T^*M)$ associated to a nonnegative quadratic Hamiltonian $H_t(\eta)=-\frac 12\langle Z_t^*J\eta,H_t^{-1}Z_t^*J\eta\rangle$.

In what follows, we identify the space $T^*_{q_0}M$ with its tangent $T_{\tilde\lambda_0}(T^*_{q_0}M)\subset
T_{\tilde\lambda_0}(T^*M)$
and introduce linear maps
$Q_t:T^*_{q_0}M\to T_{q_0}M,\ 0\le t\le 1$, by the formula: $Q_t(\eta(0))=\pi_*\eta(t)$, where $\eta(\tau),\ 0\le\tau\le t$, is a solution of equation (9) and $\eta(0)\in T^*_{q_0}M$.

\smallskip
{\bf Remark 4.} According to our assumptions, $\tilde u(\tau)$ is a strict local maximum of the function
$u\mapsto h^\tau_u(\tilde\lambda_\tau)$. Assume that this is a global maximum and moreover
$h^\tau(\lambda)=\max\limits_{u\in U}h^\tau_u(\lambda)$ is smooth with respect to $\lambda\in T^*M$.
Then $\tilde\lambda_\tau,\ 0\le\tau\le 1,$ is a solution of the Hamiltonian system
$\dot\lambda=\vec h^\tau(\lambda)$. We define the {\it exponential map} $\mathcal E^t_q:T^*_qM\to M$
by the formula $\mathcal E^t_q(\lambda_0)=\pi(\lambda_t)$, where
$\dot\lambda_\tau=\vec h^\tau(\lambda_\tau),\ 0\le\tau\le t$. It is not hard to show that
$Q_t=(\tilde P_*^t)^{-1}D_{\tilde\lambda_0}\mathcal E^t_{q_0}$.

\smallskip
\begin{theorem} Under conditions of Theorem~1, the following identities are valid:
$$
\det(I+K)=\det(Q_1\Gamma^{-1}_1),
$$
$$
\mathrm{tr}K=
\mathrm{tr}\left(\int_0^1\int_0^tX_tH^{-1}_tZ^*_tJZ_\tau H^{-1}_\tau X^*_\tau\,d\tau dt\Gamma_1^{-1}\right).
$$
\end{theorem}

Let us apply this theorem to Example~1. In our coordinates, $\Gamma_1=I$. We have to find matrix $Q_1$. It is convenient to use complex notations: $\eta=(p;q)\in\mathbb C\times\mathbb C,\ p=ip_1+p_2,\ q=iq_1+q_2$. System (8) has the form: $\dot p=irp-\nu^2q,\quad \dot q=irq+p.$
We have to find $q(t)$ under conditions $p(0)=1,\ q(0)=0$; then $Q_1$ is just complex number $q(1)$ treated as a
$2\times 2$-real matrix, $\det Q_1=|q(1)|^2$. A simple calculation gives: $q(t)=\frac{\sin\nu t}\nu e^{i\nu t}$.
Keeping in mind that all eigenvalues of the operator $K$ have multiplicity 2, we obtain the square of a classical Euler identity:$
\prod\limits_{n=1}^{+\infty}\left(1-\left(\frac r{\pi n}\right)^2\right)^2=\left(\frac{\sin r}r\right)^2.$

Now we consider one more very simple example, a harmonic oscillator.

\smallskip
{\bf Example 2.} $M=U=\mathbb R,\ f_u(q)=u,\ \ell(q,u)=\frac 12(u^2-r q^2),\ q_0=0,\ \tilde u(t)=0$. Then
$h_u(p,q)=pu-\frac 12(u^2-rq^2),\ h_{\tilde u(t)}=\frac r2q^2$,
$$
\tilde\Phi_t(p,q)=\begin{pmatrix}p-trq\\q\end{pmatrix},\
g^t_u=(p-trq)u-\frac 12 u^2,\quad
H_t=-1,\ Z_t=\begin{pmatrix}tr\\1\end{pmatrix},\ X_t=1.
$$
Operator $K$ has a form:
$$
Kv(t)=r\int_0^t(t-\tau)v(\tau)\,d\tau-r\int\limits_0^1\int\limits_0^t(t-\tau)v(\tau)\,d\tau dt.
$$
The eigenfunctions of this operator have a form $t\mapsto c\cos(\pi nt),\ c\in\mathbb R,\ n=1,2,\ldots$, where the eigenfunction $c\cos(\pi nt)$ corresponds to the eigenvalue $-\frac r{(\pi n)^2}$. Moreover, $Q_1=\frac{\sin\sqrt{r}}{\sqrt{r}}$ if $r>0$ and $Q_1=\frac{\mathrm{sh}\sqrt{|r|}}{\sqrt{|r|}}$ if $r<0$.
The determinant formula from Theorem~2 coincides with the Euler identity:
$\prod\limits_{n=1}^\infty\left(1-\frac r{(\pi n)^2}\right)=\frac{\sin\sqrt r}{\sqrt r}$ or its hyperbolic version. The trace formula gives another famous Euler observation:
$\sum\limits_{n=1}^\infty \frac r{(\pi n)^2}=\frac r6$.

Harmonic oscillator is a special case of a more general example where the spectrum is still explicitly computed.

\smallskip
{\bf Example 3.} $M=U=\mathbb R^m,\ f_u(q)=Aq+u,\ \ell(q,u)=\frac 12(|u|^2-\langle q,Rq\rangle),\ q_0=0,\ \tilde u(t)=0$, where $A$ is a $m\times m$-matrix and $R$ is a symmetric $m\times m$-matrix. Then
$h_u(p,q)=\langle p,Aq+u\rangle-\frac 12(|u|^2-\langle q,Rq\rangle),\ h_{\tilde u(t)}=\langle p,Aq\rangle+\frac 12\langle q,Rq\rangle$,
$$
\tilde\Phi_t(p,q)=\begin{pmatrix}e^{-tA^*}p-\int_0^te^{(\tau-t)A^*}Re^{\tau A}\,d\tau\\e^{tA}q\end{pmatrix},
$$
$$
g^t_u=\left\langle e^{-tA^*}p-\int_0^te^{(\tau-t)A^*}Re^{\tau A}\,d\tau,u\right\rangle -\frac 12 |u|^2,
$$
$$
H_t=-I,\quad Z_t=\begin{pmatrix}\int_0^te^{\tau A^*}Re^{\tau A}\,d\tau\\e^{-tA}\end{pmatrix},\quad X_t=e^{-tA}.
$$
Then $\mathcal V=\{v\in L^2([0,1];\mathbb R^2):\int_0^1e^{-tA}v(t)\,dt=0\}$,
$$
<Kv|v>=\int_0^1\left\langle \iint\limits_{0\le\tau_1\le\tau_2\le t}e^{(\tau_2-t)A^*}Re^{(\tau_2-\tau_1)A}v(\tau_1)\,d\tau_2d\tau_1,v(t)\right\rangle\,dt.
$$

A vector function $v(\cdot)\in\mathcal V$ is an eigenfunction of the operator $K$ with an eigenvalue $\alpha\in\mathbb R$ if and only if
$$
\alpha v(t)=\iint\limits_{0\le\tau_1\le\tau_2\le t}e^{(\tau_2-t)A^*}Re^{(\tau_2-\tau_1)A}v(\tau_1)\,d\tau_2d\tau_1+
e^{-tA^*}c,\quad t\in[0,1], \eqno (10)
$$
for some constant vector $c\in\mathbb R^m$. We denote: $y(t)=\int_0^te^{(t-\tau)A}v(\tau)\,d\tau$; then we differentiate twice equation (10) and obtain that this equation is equivalent to the boundary values problem
$y(0)=y(1)=0$ for the ordinary differential equation
$$
\alpha\ddot y=\alpha(A-A^*)\dot y+(\alpha A^*A+R)y.
$$

From now on we assume that the matrix $A$ is symmetric, $A^*=A$; then nonzero eigenvalues and corresponding eigenfunctions are nontrivial solution of the boundary value problem:
$$
\ddot y=(A^2+\frac 1\alpha R)y,\quad y(0)=y(1)=0.
$$
We obtain that $\alpha$ is a nonzero eigenvalue of $K$ if and only if there exists a positive integer $n$ such that $-(\pi n)^2$ is an eigenvalue of the matrix $A^2+\frac 1\alpha R$. Moreover, the multiplicities of the eigenvalue
$\alpha$ of $K$ and of the eigenvalue $-(\pi n)^2$ of $A^2+\frac 1\alpha R$ are equal. In other words,
$\det K=\prod\limits_{n=1}^\infty\prod\limits_i(1-s_{ni})$ where $s_{ni}$ are real roots (counted with their multiplicities) of the following polynomial equation with unknown $s$:
$$
\det(s(A^2+(\pi n)^2I)-R)=0.
$$
Actually, all roots of this polynomial are real since $A^2+(\pi n)^2I$ is a sign-definite symmetric matrix.
Let $\phi(s)=\det(s(A^2+(\pi n)^2I)-R),\ \phi(s)=\sum\limits_{i=1}^mc_{ni}s^i$; then $c_{nm}=\det(A^2+(\pi n)^2I)$ and
$\prod\limits_{n=1}^\infty\prod\limits_i(1-s_{ni})=\frac{\phi(1)}{c_{nm}}$. We obtain that
$$
\prod\limits_{i=1}^m(1-s_{ni})=\frac 1{c_{nm}}\det((A^2+(\pi n)^2I)-R)=
\det(I-R(A^2+(\pi n)^2I)^{-1}).
$$
Hence $\det K=\prod\limits_{n=1}^\infty\det(I-R(A^2+(\pi n)^2I)^{-1})$.\\
Similarly, $\mathrm{tr}K=-\sum\limits_{n=1}^\infty\mathrm{tr}(R(A^2+(\pi n)^2I)^{-1})$.

Now we have to compute the right-hand sides of the identities of Theorem~2. First of all,
$\Gamma_1=\int_0^1e^{-2tA}\,dt$. Then we consider Jacobi system (9). We have $T^*M=\mathbb R^m\times\mathbb R^m=\{(p,q): p,q\in\mathbb R^m\}$; then $Q_tp_0=q_t$ where $t\mapsto(p_t,q_t)$ is the solution of the Jacobi system with the initial value $q_0=0$. We set $y(t)=e^{tA}q_t$, differentiate in virtue of the Jacobi system and obtain the equation: $\ddot y=(A^2-R)y$. Moreover, $\dot y(0)=p_0$. Hence $Q_t=(R-A^2)^{-1/2}\sin((A^2-R)^{1/2}t)$. This formula is valid also for a sign-indefinite matrix $R-A^2$ if we properly interpret the square root or simply make the computations in coordinates where the matrix $A^2-R$ is diagonal. Finally, putting together all the formulas, we obtain the following generalization of the Euler identities:

\begin{prop} Let $A,R$ be symmetric matrices. Then:
$$
\prod\limits_{n=1}^\infty\det\left(I-R(A^2+(\pi n)^2I)^{-1}\right)=\frac{2\det\left(\sin\sqrt{R-A^2}\right)}{\det\left(\sqrt{R-A^2}\int_{-1}^1e^{tA}\,dt\right)},
$$
$$
\sum\limits_{n=1}^\infty\mathrm{tr}(R(A^2+(\pi n)^2I)^{-1})=
$$
$$
\mathrm{tr}\bigl(\iiint\limits_{0\le\tau_1\le\tau_2\le t\le 1}e^{(\tau_2-2t)A}Re^{(\tau_2-2\tau_1)}\,d\tau_2d\tau_1dt(\int_0^1e^{-2tA}\,dt)^{-1}\bigr).
$$
\end{prop}
The right-hand side of the determinant formula has an obvious meaning also in the case of a degenerate $R-A^2$.
If $m=1,\ A=a,\ R=a^2+b^2$, we get:
$$
\prod\limits_{n=1}^\infty\left(1-\frac{a^2+b^2}{a^2+(\pi n)^2}\right)=
\frac{a\sin b}{b\,\mathrm{sh}\,a},
$$
an interpolation between the classical Euler identity and its hyperbolic version.
The trace identity is essentially simplified if the matrices $R$ and $A$ commute. In the commutative case we obtain:
$$
\sum\limits_{n=1}^\infty\mathrm{tr}(R(A^2+(\pi n)^2I)^{-1})=
\frac 12\mathrm{tr}\left(R(A\,\mathrm{cth}A-I)A^{-2}\right).
$$

\section{Proof of Theorem 1}

We start with some definitions and notations. A compact quadratic form $b$ on the Hilbert space $\mathcal V$ is a form defined by a compact symmetric operator $B,\ b(v)=<Bv|v>,\ v\in\mathcal V$. The spectrum of $b$ and its positive and negative parts are, by the definition, those of $B$, i.\,e.
$\mathrm{Sp}_\pm(b)=\mathrm{Sp}_\pm(B)$. Recall that the eigenvalues of $B$ are just critical values of the restriction of $b$ to the unit sphere in $\mathcal V$.

Let $\mathcal V_0\subset\mathcal V$ be a Hilbert subspace, then  the form $b|_{\mathcal V_0}$ is defined by the composition of $B|_{\mathcal V_0}$ and the orthogonal projection of $\mathcal V$ on $\mathcal V_0$. The relation between $\mathrm{Sp}(b)$ and $\mathrm{Sp}(b|_{\mathcal V_0})$ is ruled by the classical Rayleigh--Courant minimax principle (see \cite{CoHi,Ka}).

Assume that $\mathrm{Sp}_\pm(b)$ is infinite; then a natural ordering $\beta_n,\ n\in\mathbb Z_\pm$, of $\mathrm{Sp}_\pm(b)$ is defined in the same way as the natural ordering of $\mathrm{Sp}_\pm(K)$ (see (7)).

\begin{definition} We say that $b$ has the spectrum of capacity $\varsigma>0$ with the remainder of order $\nu>1$ if
$\mathrm{Sp}_+(b)$ and $\mathrm{Sp}_-(b)$ are both infinite and
$$
\beta_n=\frac \varsigma n+O(n^{-\nu})\quad \mathrm{as}\quad n\to\pm\infty. \eqno (11)
$$
We say that $b$ has the spectrum of zero capacity with the remainder of order $\nu$ if either $\mathrm{Sp}_\pm(b)$
is finite or $\beta_n=O(n^{-\nu})$ as $n\to\pm\infty$.
\end{definition}

Let $b_i$ be a quadratic form on the Hilbert space $\mathcal V_i,\ i=1,2$; then $b_1\oplus b_2$ is a quadratic
form on $\mathcal V_1\oplus\mathcal V_2$, $\mathrm{Sp}(b_1\oplus b_2)=\mathrm{Sp}(b_1)\cup\mathrm{Sp}(b_2)$
and the multiplicities of common eigenvalues are added.

\begin{prop}
\begin{description}
\item[(i)] If $b$ has the spectrum of capacity $\varsigma\ge 0$, then $sb$ has the spectrum of capacity $s\varsigma$ with the remainder of the same order as $b$, for any $s\in\mathbb R$.
\item[(ii)] If $b_1,b_2$ have the spectra of capacities $\varsigma_1,\varsigma_2$ with the remainders of equal orders, then $b_1\oplus b_2$ has the spectrum of capacity $\varsigma_1+\varsigma_2$ and the remainder of the same order as $b_1,b_2$.
\item[(iii)] Let $\mathcal V_0$ be a Hilbert subspace of the Hilbert space $\mathcal V$ and
$\dim(\mathcal V/\mathcal V_0)<\infty$. Assume that one of two forms  $b$ or $b|_{\mathcal V_0}$ has the spectrum of capacity $\varsigma\ge 0$ with a remainder of order $\nu\le 2$. Then the second form has the spectrum of the
the same capacity $\varsigma$ with a remainder of the same order $\nu$.
\item[(iv)] Let the forms $b$ and $\hat b$ be defined on the same Hilbert space $\mathcal V$, where $b$ has the spectrum of capacity $\varsigma$ and $\hat b$ has the spectrum of zero capacity, both with the reminder term of order $\nu$. Then the form $b+\hat b$ has the spectrum of capacity $\varsigma$ with the reminder term of order $\frac{2\nu+1}{\nu+1}$.
\end{description}
\end{prop}

{\bf Proof.} Statement (i) is obvious.  To prove (ii)
we re-write asymptotic relation (11) in a more convenient form. An equivalent relation for positive $n$
reads:
$$
\#\left\{k\in\mathbb Z: 0<\frac 1{\beta_k}<n\right\}=\varsigma n+O(n^{2-\nu}),\quad \mathrm{as}\quad n\to\infty
$$
and similarly for negative $n$. Statement (ii) follows immediately.

Statement (iii) follows from the Rayleigh--Courant minimax principle for the eigenvalues and the
relation: $\left|\frac{\varsigma}n-\frac{\varsigma}{n+j}\right|=O(\frac 1{n^2})$ as $|n|\to\infty$ for any fixed $j$.

To prove (iv) we use the Weyl inequality for the eigenvalues of the sum of two forms. Weyl inequality is a straightforward
corollary of the minimax principle, it claims that the positive eigenvalue number $i+j-1$ in the natural ordering of the sum of two forms does not exceed the sum of the eigenvalue number $i$ of the first summand and the eigenvalue number $j$ of the second summand. Of course, we may equally works with naturally ordered negative eigenvalues simply changing the signs of the forms.

In our case, to have both sides estimates we first present $b+\hat b$ as
the sum of $b$ and $\hat b$ and then present $b$ as the sum of $b+\hat b$ and $-\hat b$. In the first case we apply the Weyl inequality with
$i=n-[n^\delta],\ j=[n^\delta]$ for some $\delta\in(0,1)$, and in the second case we take $i=n,\ j=[n^\delta]$.
The best result is obtained for $\delta=\frac 1{\nu+1}. \qquad\square$

\smallskip
We have to prove that the spectrum of operator $K$ (see (6), (5)) has capacity $\frac 1{\pi}\int_0^1\bar\zeta_t\,dt$ with the remainder of order $\frac 53$. First, we may identify $T_{\tilde u(t)}$
with $\mathbb R^m$ and assume that $H_t=-I$. Indeed, if we trivialize the vector bundle $\tilde u^*(TU)$ over the
segment $[0,1]$, then $H_t$ becomes a negative definite symmetric matrix. Then we substitute $v$ by $(-H_t)^{\frac 12}v$ and $Z_t$ by $Z_t(-H_t)^{-\frac 12}$.

Let $0=t_0<t_1<\cdots<t_l<t_{l+1}=1$ be a subdivision of the segment $[0,1]$. The subspace
$$
\{v\in\mathcal V: \int_{t_i}^{t_{i+1}}X_tv(t)\,dt=0,\ i=0,1,\ldots,l\}=\bigoplus\limits_{i=0}^l\mathcal V_i,\
\mathcal V_i\subset L^2([0,1];\mathbb R^m)
$$
has a finite codimension in $\mathcal V$. The quadratic form
$$
<Kv|v>=\int_0^1\left\langle JZ_tv(t),\int_0^tZ_\tau v(\tau)\,d\tau\right\rangle\,dt \eqno (12)
$$
restricted to this subspace turns into the direct sum of the forms
$$
<K_iv|v>=\int_{t_i}^{t_{i+1}}\left\langle JZ_tv(t),\int_{t_i}^tZ_\tau v(\tau)\,d\tau\right\rangle\,dt,\quad
v_i\in\mathcal V_i,\ i=0,1\ldots l.
$$
Indeed, the relations $\int_{t_i}^{t_{i+1}}X_tv(t)\,dt=0$ imply that
$$
\left\langle\int_{t_i}^{t_{i+1}}JZ_tv(t)\,dt,\int_{t_j}^{t_{j+1}}Z_tv(t)\,dt\right\rangle=0.
$$
According to Proposition~3 it is sufficient to prove our theorem for the operators $K_i,\ i=0,1\ldots,l$.
In particular, we may substitute the piecewise analyticity assumption in the statement of the theorem by the analyticity one.

Moreover, under the analyticity condition we may assume that
$$
Z_t^*JZ_t=\bigoplus\limits_{j=1}^k\begin{pmatrix}0& -\zeta_j(t)\\ \zeta_j(t)& 0\end{pmatrix},
$$
where $0\le 2k\le m$ and $\zeta_j(t)$ are not identical zero. Indeed, according to the Rayleigh theorem (see \cite{Ka}), there exists an analytically depending on $t$ orthonormal  basis in which our anti-symmetric matrix takes a desired form.

The functions $\zeta_j(t),\ j=1,\ldots,k$, are analytic and may have only isolated zeros. Hence we may take a finer subdivision of $[0,1]$ in such a way that $\zeta_j(t),\ j=1,\ldots,k$, do not change sign on the segments $[t_i,t_{i+1}]$. Moreover, a simple change of the basis of $\mathbb R^m$ if necessary allows us to assume that
$\zeta_j(t)\ge 0,\ t\in[t_i,t_{i+1}]$. Actually, to simplify notations a little bit, we may simply assume that
$\zeta_j(t)\ge 0,\ 0\le t\le 1,\ j=1,\ldots,k$. In this case $\bar\zeta(t)=\sum\limits_{j=1}^k\zeta_j(t)$.

Let us study quadratic form (12) on the space $\{v\in L^2([0,1];\mathrm R^m):\int_0^1v(t)\,dt=0\}$.
Recall that we are allowed by Proposition~3 to work on any subspace of $L^2([0,1];\mathrm R^m)$ of a finite
codimension. We set $w(t)=\int_0^tv(\tau)\,d\tau$; a double integration by parts gives:
$$
\int_0^1\left\langle JZ_tv(t),\int_0^tZ_\tau v(\tau)\,d\tau\right\rangle\,dt=
\int_0^1\left\langle JZ_tv(t),Z_tw(t)\right\rangle\,dt+
$$
$$
\int_0^1\left\langle JZ_tw(t),\dot Z_tw(t)\right\rangle\,dt+\int_0^1\left\langle J\dot Z_tw(t),\int_0^t\dot Z_\tau w(\tau)\,d\tau\right\rangle\,dt.
$$
Moreover, we have:
$$
\left|\int_0^1\left\langle JZ_tw(t),\dot Z_tw(t)\right\rangle\,dt+\int_0^1\left\langle J\dot Z_tw(t),\int_0^t\dot Z_\tau w(\tau)\,d\tau\right\rangle\,dt\right|\le c\int_0^1|w(t)|^2\,dt
$$
for some constant $c$.

Let $\bar\lambda_n$ and $\underline{\lambda}_n,\ n\in\mathbb Z\setminus\{0\}$, be naturally ordered non zero eigenvalues of the quadratic forms
$$
\int_0^1\left\langle JZ_tv(t),Z_tw(t)\right\rangle\,dt+c\int_0^1|w(t)|^2\,dt
$$
and
$$
\int_0^1\left\langle JZ_tv(t),Z_tw(t)\right\rangle\,dt-c\int_0^1|w(t)|^2\,dt
$$
correspondently. The minimax principle for the eigenvalues implies that $\underline{\lambda}_n\le\alpha_n\le\bar\lambda_n,\ n\in\mathbb Z\setminus\{0\}$. Moreover, the form
$$
\int_0^1\left\langle JZ_tv(t),Z_tw(t)\right\rangle\,dt\pm c\int_0^1|w(t)|^2\,dt
$$
splits in the direct sum of the forms
$$
\int_0^1\zeta_j(t)\langle Jv_j(t),w_j(t)\rangle\,dt\pm c\int_0^1|w_j(t)|^2\,dt,\quad j=1,\ldots,k, \eqno (13)
$$
where $v_j(t)\in\mathbb R^2,\ w_j(t)=\int_0^tv_j(t),\ w_j(1)=0$, $J=\bigl(\begin{smallmatrix}0 & -1\\1 & 0\end{smallmatrix}\bigr)$ and we simply ignore the identically vanishing part if $2k<m$.

It remains to estimate the spectrum of the forms (13). To do that, we study the spectrum of the forms
$$
\int_0^1\zeta_j(t)\langle Jv_j(t),w_j(t)\rangle\,dt,\qquad \int_0^1|w_j(t)|^2\,dt, \eqno (14)
$$
and apply the statement (iv) of Proposition~3 to their linear combinations. To simplify calculations, we identify
$\mathbb R^2$ with $\mathbb C$ as follows: $(v_j^1,v_j^2)=v_j^1+iv_j^2$, then $J$ is the multiplication on the
imaginary unit $i$. A complex-valued function $v_j(\cdot)$ is an eigenfunction  with an eigenvalue $\lambda$ for the first of two forms (14) if and only if it is a critical point of the functional
$$
\int_0^1\langle\zeta(t)iv_j(t),w_j(t)\rangle-\lambda\langle v_j(t),v_j(t)\rangle\,dt\quad w_j(0)=w_j(1)=0.
$$
The Euler--Lagrange equation for this functional reads:
$$
2\lambda\ddot w_j+(\zeta_j iw_j)^\cdot+\zeta_j i\dot w=0,\quad w_j(0)=w_j(1)=0.
$$
We set $x(t)=e^{\frac{i}{2\lambda}\int_0^t\zeta_j(\tau)\,d\tau} w_j(t)$, plug-in this expression in the equation and
arrive to a standard Sturm--Liouville problem:
$$
\ddot x+\left(\frac{\zeta_j}{2\lambda}\right)^2x=0,\quad x(0)=x(1)=0.
$$
We see that the spectrum is double (recall that $x\in\mathbb C$) and is symmetric with respect to the origin.
The asymptotics of the spectrum for the Sturm--Liouville problem is well-known (see \cite{CoHi}). We obtain:
$$
\lambda_n=\frac 1{\pi n}\int_0^1\zeta_j(t)\,dt+O\left(\frac 1{n^2}\right),\quad n\to\pm\infty.
$$
In other words, this spectrum has capacity $\frac 1{\pi}\int_0^1\zeta_j(t)\,dt$ and the remainder term of order 2.
The Euler-Lagrange equation for the eigenfunctions problem for the second of two forms (14) reads:
$\lambda\ddot w_j+w_j=0,\ w_j(0)=w_j(1)=0$; hence $\lambda_n=\frac 1{(\pi n)^2}=O(\frac 1{n^2})$. The spectrum of this form has zero capacity with the remainder of order 2. $\qquad\square$

\section{Proof of Theorem 2}

We again assume that $H_t=-I$ using the same preliminary change of variables as in the proof of Theorem~1, if necessary: substitute $v$ by $(-H)^{\frac 12}v$ and $Z_t$ by $Z_t(-H)^{-\frac 12}$. Moreover, we fix some coordinates in a neighborhood of $q_0\in M$ and use induced coordinates in $T^*M$ so that $T_{\tilde\lambda_0}(T^*M)$ is identified
with $\mathbb R^d\times\mathbb R^d=\{(p,q):p,q\in\mathbb R^d\}$, where $d=\dim M,\ \pi_*(p,q)=q$, and
$$
\sigma_{\tilde\lambda_0}\left((p,q),(p',q')\right)=\left\langle J\begin{pmatrix}p\\q\end{pmatrix},\begin{pmatrix}p'\\q'\end{pmatrix}\right\rangle,\quad J=\begin{pmatrix}0& -I\\I& 0\end{pmatrix}.
$$

The map $Z_t:\mathbb R^m\to\mathbb R^d\times\mathbb R^d$ has a form $Z_t:v\mapsto (Y_tv,X_tv)$, where $Y_t$ and $X_t$ are $d\times m$-matrices. Let $s\in\mathbb C$; we define $Z^s_t:\mathbb C^d\times C^d$ by the formula
$Z_t^sv=(sY_tv,X_tv)$. Now consider the complexified Jacobi equation:
$$
\dot\eta=Z^s_t{Z^s_t}^*J\eta,\quad \eta\in\mathbb C^d\times\mathbb C^d, \eqno(15),
$$
where the transposition ``\,*\,'' corresponds to the complex inner product and not to the Hermitian one! The matrix form of this equation is as follows:
$$
\begin{pmatrix}\dot p\\ \dot q\end{pmatrix}= \begin{pmatrix}sY_t\\X_t\end{pmatrix}\bigl(sY_t^*,X_t^*\bigr)\begin{pmatrix}-q\\p\end{pmatrix}. \eqno (16)
$$

\begin{prop}
Let $\Phi_t^s:\mathbb C^d\times\mathbb C^d\to \mathbb C^d\times\mathbb C^d,\ 0\le t\le 1$ be the fundamental matrix
of the complexified Jacobi equation, i.\,e.
$$
\frac d{dt}\Phi^s_t=Z^s_t{Z^s_t}^*J\Phi^s_t,\quad \Phi^s_0=I;
$$
then $s\mapsto\Phi_1^s,\ s\in\mathbb C,$ is an entire matrix-valued function and
$$
\|\Phi_1^s\|\le ce^{c|s|}
$$
for some constant $c$
\end{prop}
{\bf Proof.} The Volterra series for our system has the form:
$$
\Phi^s_1=I+\sum\limits_{n=1}^\infty\int\limits_0^1\int\limits_0^{t_1}\cdots\int\limits_0^{t_{n-1}}
Z^s_{t_1}{Z^s_{t_1}}^*J\cdots Z^s_{t_n}{Z^s_{t_n}}^*J\,dt_n\ldots dt_1 .
$$
The $n$-th term of the series is a polynomial matrix of $s$ whose norm is bounded by
$\frac 1{n!}\left(\int_0^1\|Z_t^s{Z_t^s}^*\|\,dt\right)^n$ for any $s$ and $n$. The sum of the series is an entire function by the Weierstrass theorem. Moreover, we obtain the estimate:
$\|\Phi_t^s\|\le e^{\int_0^1\|Z_t^s{Z_t^s}^*\|\,dt}$ but it is worse than one we need.

To obtain a better estimate we make a change of variables in the matrix form of the complexified Jacobi equation (see (16)). We set $\eta^s=(p,sq)$; then $\dot\eta^s=sZ_tZ_t^*J\eta^s$ and we obtain:
$|\eta^s(1)|\le e^{|s|\int_0^1\|Z_tZ_t^*\|\,dt}|\eta^s(0)|$
that gives us the desired estimate of $\|\Phi^s_1\|$ for $s$ separated from zero.$\qquad \square$

\smallskip
We define linear maps $Q^s_t:\mathbb C^n\to\mathbb C^n$ by the formula: $Q^s_tp(0)=q(t)$, where
$t\mapsto\left(\begin{smallmatrix}p(t)\\q(t)\end{smallmatrix}\right)$ is the solution of system (16) with the initial
condition $q(0)=0$. Note that the matrix $Q^s_t$ is real if $s\in\mathbb R$; moreover, $Q^1_t\bigr|_{\mathbb R^d}=Q_t$ (linear maps $Q_t$ were defined just before the statement of Theorem~2). Proposition~4 implies that
$\|Q_1^s\|\le ce^{c|s|}$.

\begin{prop} Let $s\in\mathbb R\setminus\{0\}$; the matrix $Q^s_1$ is degenerate if and only if $-\frac 1s\in\mathrm{Sp}(K)$; moreover,
$\dim\ker Q_1^s$ equals the multiplicity of the eigenvalue $-\frac 1s$.
\end{prop}
{\bf Proof.} Let $v\in\mathcal V$; it is easy to see that
$$
\left\langle JZ_\tau^sv(\tau),Z_t^sv(t)\right\rangle=s\left\langle JZ_\tau v(\tau),Z_tv(t)\right\rangle.
$$
Hence
$$
<(I+sK)v|v>=\int_0^1\langle v(t),v(t)\rangle\,dt-\int_0^1\left\langle J\int_0^tZ_\tau^sv(\tau)\,d\tau,Z_t^sv(t)\right\rangle\,dt.
$$
We keep symbol $\mathcal V$ for the complexification of $\mathcal V$, i.\,e. $\mathcal V=\mathbb C\otimes\mathcal V$. A vector-function $v\in\mathcal V$ is an eigenvector of $K$ with the eigenvalue $-\frac 1s$ if and only if
$$
\int_0^1\left\langle v(t)-{Z_t^s}^*J\int_0^tZ_t^sv(\tau)\,d\tau,w(t)\right\rangle\,dt=0,\quad \forall\,w\in\mathcal V, \eqno(17)
$$
where $\mathcal V=\left\{w\in L_2([0,1];\mathbb C^d) : \int_0^1X_tw(t)\,dt=0\right\}$.

Equation $\int_0^1X_tw(t)\,dt=0$ can be rewritten as follows: $\int_0^1\langle J\nu,Z_tw(t)\rangle\,dt=0$ for any
$\nu\in\ker\pi_*$. In other words,
$\mathcal V^\perp=\{t\mapsto {Z_t^s}^*J\nu : \nu\in\ker\pi_*\}$.
Hence relation (17) is equivalent to the existence of $\nu\in\ker\pi_*$ such that
$$
v(t)={Z_t^s}^*J\left(\int_0^tZ_\tau^sv(\tau)\,d\tau+\nu\right).
$$
Moreover, the vector $\nu$ is unique for given vector-function $v$ since the Volterra equation
$v(t)={Z_t^s}^*J\int_0^tZ_\tau^sv(\tau)\,d\tau$ has only zero solution. We set
$\eta(t)=\int_0^tZ^s_\tau v(\tau)\,d\tau+\nu$ and obtain that a vector-function $v$ satisfies relation (17) if and
only if
$$
\dot\eta=Z_t^s{Z_t^s}^*J\eta,\quad 0\le t\le 1,\quad \eta(0),\eta(1)\in\ker\pi_*. \eqno (18)
$$

It follows that $\dim\ker Q^s_1$ is equal to the multiplicity of the eigenvalue $-\frac 1s$ of the operator $K$ plus the dimension of the space of constant solutions of equation (18) that belong to $\ker\pi_*$. Let
$(\nu,0)\in\mathbb C^d\times\mathbb C^d$ be such a solution. We plug-in it in (18) and get
$X_tX_t^*\nu=0,\ \forall\,t\in[0,1]$. Hence $\nu=0.\qquad \square$

\begin{corollary} The equation $\det Q^s_1=0,\ s\in\mathbb C$, has only real roots.
\end{corollary}
Indeed, the operator $K$ is symmetric and has only real eigenvalues.$\qquad\square$

\begin{prop} Let $s_0\in\mathbb R\setminus\{0\}$ and $\det Q_1^{s_0}=0$. Then the multiplicity of the root $s_0$ of the equation
$\det Q_1^s=0$ is equal to $\dim\ker Q_1^{s_0}$.
\end{prop}
{\bf Proof.} We may assume that $s\in\mathbb R$ and work in the real setting. We denote by
$\eta^s(t)=(p^s(t),q^s(t)),\ t\in[0,1],$ those solutions of the Hamiltonian system
$$
\dot\eta^s=sZ_tZ_t^*J\eta^s \eqno (17)
$$
that satisfy the initial condition $q^s(0)=0$; then $q^s(t)=sQ^s_tp^s(0)$ (c.\,f. proof of Proposition~4).

Let $\Lambda_1^s=\{(p^s(1),q^s(1)): p^s(0)\in\mathbb R^d\}$; then $\Lambda_1^s$ is a Lagrange subspace of the symplectic space $\mathbb R^d\times\mathbb R^d$ endowed with the standard symplectic structure
$\left((p_1,q_1),(p_2,q_2)\right)\mapsto\langle p_1,q_2\rangle-\langle p_2,q_1\rangle$.

We take another Lagrange subspace $\Delta$ that is transversal to $\ker\pi^*_0$ and $\Lambda_1^*$. Then
$\Delta=\{(Aq,q):q\in\mathbb R^d\}$ where $A$ is a symmetric matrix. We make a symplectic change of variables
$(p,q)\mapsto(p',q)$ by putting $p'=p-Aq$. In new variables, $\Lambda_1^s$ is transversal to the ``horizontal''
subspace defined by the equation $p'=0$ for all $s$ close to $s_0$. Hence
$$
\Lambda_1^s=\{(p',R(s)p'):p'\in\mathbb R^d\},
$$
where $R(s)$ is a symmetric matrix. According to our construction, $sQ_1^s$ is the product of an analytic with respect to $s$ nondegenerate matrix and $R(s)$. It follows that $\dim\ker R(s)=\dim\ker Q^s_1$ and the multiplicities of the root $s_0$ of the equations $\det R(s)=0$ and $\det Q_1^s=0$ are equal. We denote
this multiplicity by $\mu(s_0)$.

Standard perturbation theory for the eigenvalues of symmetric operators (see, for instance, \cite{Ka}) implies that: (i) $\mu(s_0)\ge\dim\ker R(s_0)$; (ii) if $\det\frac{dR}{ds}(s_0)\ne 0$, then $\mu(s_0)=\dim\ker R(s_0)$.

Now we have to compute $\frac{dR}{ds}$. To do that, we re-write equation (17) in coordinates $(p',q)$. In new coordinates, this Hamiltonian equation has the same structure as (17) with $Z_t=(Y_t,X_t)$ substituted by
$Z'_t=(Y_t+AX_t,X_t)$. From now on, we'll write $p, Z_t$ instead of $p',Z_t'$ and work directly with equation (17) in order to simplify notations. We have $\eta^s(1)=(p^s(1),R(s)q^s(1))$; the symmetricity of the matrix
$R(s)$ implies that
$$
\left\langle J\eta^s(1),\frac{\partial\eta^s}{\partial s}(1)\right\rangle=\left\langle p^s(1),\Bigl(\frac{dR}{ds}\Bigr) p^s(1)\right\rangle.
$$
Proposition~6 is now reduced to the following lemma:

\begin{lemma} If $\left\langle J\eta^s(1),\frac{\partial\eta^s}{\partial s}(1)\right\rangle=0$, then $\eta^s(1)=0$.
\end{lemma}
{\bf Proof.} We have: $
\frac{\partial\dot\eta^s}{\partial s}=sZ_tZ_t^*J\frac{\partial\eta^s}{\partial s}+Z_tZ_t^*J\eta^s.$
Hence
$$
\frac\partial{\partial t}\left\langle J\eta^s(t),\frac{\partial\eta^s}{\partial s}(t)\right\rangle=
s\left\langle JZ_tZ^*_tJ\eta^s(t),\frac{\partial\eta^s}{\partial s}(t)\right\rangle+
s\left\langle J\eta^s(t),Z_tZ^*_tJ\frac{\partial\eta^s}{\partial s}(t)\right\rangle
$$
$$
+\left\langle J\eta^s(t),Z_tZ^*_tJ\eta^s\right\rangle.
$$
First two terms in the right hand side of the last identity cancel because the matrix $J$ is anti-symmetric and the matrix $JZ_tZ_t^*J$ is symmetric. We obtain:
$$
\frac\partial{\partial t}\left\langle J\eta^s(t),\frac{\partial\eta^s}{\partial s}(t)\right\rangle=
\left\langle Z_t^*J\eta^s(t),Z_t^*J\eta^s(t)\right\rangle.
$$
Moreover, $\left\langle J\eta^s(0),\frac{\partial\eta^s}{\partial s}(0)\right\rangle=0$ since $q^s(0)=0$ for all $s$. We get:
$$
\left\langle J\eta^s(1),\frac{\partial\eta^s}{\partial s}(1)\right\rangle=
\int_0^1|Z_t^*J\eta^s(t)|^2\,dt.
$$
If $\left\langle J\eta^s(1),\frac{\partial\eta^s}{\partial s}(1)\right\rangle=0$, then
$Z^*_tJ\eta^s(t)=0,\ \forall\,t\in[0,1]$. Equation (17) now implies that $\dot\eta^s(t)\equiv 0$ and
$\eta^s(t)=(p^s(0),0)$. We plug-in this equilibrium in (17) and obtain: $X_tX_t^*p^s(0)=0,\ 0\le t\le 1$.
Hence $\Gamma_1p^s(0)=\int_0^1X_tX_t^*\,dtp^s(0)=0$. Recall that $\Gamma_1$ is a nondegenerate matrix.
$\qquad\square$

\smallskip
Proposition~4 implies that $s\mapsto\det Q_1^s,\ s\in\mathbb C$, is an entire function that satisfies the estimate
$\det Q_1^s\le ce^{c|s|}$ for some constant $c$. It follows from a classical Hadamard theorem (see \cite{Co}) that such a function has a presentation $\det Q^s_1=ae^{bs}\prod_i(1-\frac s{s_i})$, where $s_i$ are roots of the equation
$\det Q_1^s=0$ counted according to multiplicity and $a,b\in\mathbb C$.

Propositions 5 and 6 now imply that
$$
\det Q_1^s=ae^{bs}\det(I+sK),\quad s\in\mathbb C.  \eqno(18)
$$
It remains to find $a$ and $b$.
We have: $\det Q_1^0=a$. Moreover,
$$
\frac d{ds}\det Q_1^s|_{s=0}=\mathrm{tr}\left(\frac{dQ^s_1}{ds}\bigr|_{s=0}(Q^0_1)^{-1}\right)\det Q_1^0.
$$
We differentiate identity (18) and obtain:
$\mathrm{tr}\left(\frac{dQ^s_1}{ds}\bigr|_{s=0}(Q^0_1)^{-1}\right)=b +\mathrm{tr}\,K$.

Let $\Pi:L^2([0,1];\mathbb R^m)\to\mathcal V$ be the orthogonal projection, then
$\mathrm{tr}\,K=\mathrm{tr}\,(\Pi\hat K\Pi)=\mathrm{tr}\,(K\Pi)$, where
$\hat Kv(t)=-\int_0^tZ_t^*JZ_\tau v(\tau)\,d\tau$.
The subspace $\mathcal V$ is the kernel of the operator $A:v\mapsto\int_0^1X_tv(t)\,dt,\ v\in\L^2([0,1];\mathbb R^m)$. We have:\quad
$\Pi=I-A^*(AA^*)^{-1}A,\quad \Pi v(t)=v(t)-X_t^*\Gamma_1^{-1}\int_0^1X_\tau v(\tau)\,d\tau.$

The trace of a trace-class operator $\bar B$ of a form $(\bar Bv)(t)=\int_0^1B(t,\tau)v(\tau)\,d\tau$ is
computed according to the formula $\mathrm{tr}\,\bar B=\int_0^1\mathrm{tr}\,B(t,t)\,dt$. This presentation of the
trace is valid in our situation as well. Moreover, $\mathrm{tr}\,(Z_t^*JZ_t)\equiv 0$; we collect the terms and
obtain the formula:
$$
\mathrm{tr}\,K=\int_0^1\int_0^t\mathrm{tr}\,(Z_t^*JZ_\tau X_\tau^*\Gamma_1^{-1}X_t)\,d\tau dt
$$
$$
=\int_0^1\int_0^t\mathrm{tr}\,(X_tZ_t^*JZ_\tau X_\tau^*\Gamma_1^{-1})\,d\tau dt
$$
as in the statement of Theorem~2 (recall that we are working in coordinated where $H_t=-I$).

Now we compute $Q^0_1$. For $s=0$, system (18) is reduced to the system: $\dot p=0,\ \dot q=X_tX_t^*p$.
Hence $Q^0_1=\int_0^1X_tX_t^*\,dt=\Gamma_1$.

The last step is the calculation of $\frac{dQ^s_1}{ds}\bigr|_{s=0}$. System (16) reads:
$$
\dot p=sY_tX_t^*p-s^2Y_tY_t^*q
$$
$$
\dot q=X_tX_t^*p-sX_tY_t^*q.
$$
We have:
$$
q^s(t)=\int_0^tX_\tau X^*_\tau\,d\tau p_0+O(s),\quad p^s(t)=s\int_0^tY_\tau X_\tau^*\,d\tau p_0+O(s^2).
$$
Hence
$$
q^s(1)=\int_0^1X_tX_t^*p(t)^s-sX_tY_t^*q^s(t)\,dt=
$$
$$
\Gamma_1p_0+s\int_0^1\int_0^tX_t X_t^*Y_\tau X_\tau^*-X_t Y_t^*X_\tau X_\tau^*\,d\tau dtp_0+O(s^2).
$$
Moreover, $X_t X_t^*Y_\tau X_\tau^*-X_t Y_t^*X_\tau X_\tau^*=X_tZ_t^*JZ_\tau X_\tau^*$. We see that \linebreak
$\mathrm{tr}\left(\frac{dQ^s_1}{ds}\bigr|_{s=0}\Gamma_1^{-1}\right)=\mathrm{tr}\,K.$ Hence $b=0$.
$\qquad\square$

\smallskip
{\bf Remark 5.} We used only the first derivative of $Q^s_1$ at $s=0$ in the proof of Theorem~2 but system (16)
allows us to find explicit integral expressions for all higher derivatives and thus to obtain integral expressions
for all elementary symmetric functions of the eigenvalues of the operator $K$.


\begin{thebibliography}{9}

\bibitem
{AgSa} A. Agrachev, Yu. Sachkov, {\it Control Theory from the geometric viewpoint.}
Springer Verlag, 2004, 426\,p.

\bibitem{BoTr} S. Bolotin, D. Treschev, {\it Hill's formula.} Russian Math. Surveys, 2010, v.\,65, 191--257

\bibitem{Co} J. Conway, {\it Functions on one complex variable, I.} Springer Verlag, 1995, 317\,p.

\bibitem{CoHi} R. Courant, D. Hilbert, {\it Methods of mathematical physics, I.} Intersience Pub., 1953,
576\,p.

\bibitem{Ka} T. Kato, {\it Perturbation theory for linear operators.} Springer Verlag, 1980, 643\,p.

\end{thebibliography}
\end{document}